\documentclass[prl, reprint, aps]{revtex4-2}

\usepackage{amsmath,amssymb,amsfonts}
\usepackage{algorithmic}
\usepackage{graphicx}
\usepackage{textcomp}
\usepackage{xcolor}
\usepackage{natbib}

\graphicspath{{./}}

\usepackage{braket}

\usepackage[ruled, vlined, linesnumbered, boxed]{algorithm2e}

\usepackage{xspace}

\newcommand{\schr}{Schr\"{o}dinger\xspace}

\newcommand{\xv}{\mathbf{x}}

\newcommand{\bv}{\mathbf{b}}

\begin{document}

\preprint{LLNL-JRNL-841953}

\title{Optimal control from inverse scattering via single-sided focusing}
\author{Michael D. Schneider}
\email[Corresponding author:]{schneider42@llnl.gov}
\author{Caleb Miller}
\author{George F. Chapline} 
\author{Jane Pratt}
\author{Dan Merl}
\affiliation{
	Lawrence Livermore National Laboratory, 7000 East Avenue, Livermore, CA 94550, USA.
}

\date{\today}

% ---------------------------------------------------------
% Feynman showed we can use PIs to solve Shcrodinger. We're using Schrodinger 
% to solve PIs (for control).
% We have the exact correspondence between Schrodinger and control worked out for the first time

\begin{abstract}
% Feynman originally showed how to use path integrals to solve the \schr equation. Here we use the \schr equation to solve Euclidean path integral solutions of optimal control problems.

We describe an algorithm to solve Bellman optimization that replaces a sum over paths determining the optimal cost-to-go by an analytic method localized in state space. Our approach follows from the established relation between stochastic control problems in the class of linear Markov decision processes and quantum inverse scattering. We introduce a practical online computational method to solve for a potential function that informs optimal agent actions. 
% The potential function is related to a state cost function for the controlled agent and defines a dynamic trajectory for the optimal control solution.
This approach suggests that optimal control problems, including those with many degrees of freedom, can be solved with parallel computations.
% determined from a least action principle. 
\end{abstract}

\maketitle

% ---------------------------------------------------------
\paragraph{Introduction}
\label{sec:introduction}
Optimal control of noisy systems~\citep{stengel1994optimal,alekseev2013optimal} 
arises in numerous applications including robotics~\citep{abdallah1991survey}, 
production planning~\citep{ivanov2012applicability}, 
power systems~\citep{christensen2013optimal}, 
traffic management~\citep{gugat2005optimal}, 
and financial investments~\citep{bertsimas1998optimal}. 
Indeed the introduction of noise to intrinsically noiseless physical systems 
allows for efficient solution by optimal control methods.
In addition, optimal control is closely related to reinforcement learning 
(RL)~\citep{sutton2018reinforcement}, which has seen a surge of research 
activity and applications in the past decade, spurred in part by computational 
advances and new more scalable solution algorithms. 
New methods to solve optimal control problems might thus have wide ranging 
impacts in both traditional control applications and RL-based machine learning 
problems. 

In the dynamic programming approach to solving optimal control problems~\citep{bertsekas2012dynamic},
one seeks to optimize the expected accumulated `cost' along state space 
trajectories that reach a target state. The immediate costs incurred 
can include terms for both state and action 
costs~\citep[{\it e.g.},][]{berkovitz2013optimal}. The choice of cost functions has often been determined by mathematical 
tractability or by heuristic assumptions about the desired performance of the
control solutions.
In robotics applications with high-dimensional state spaces, cost function 
specification has been adapted to the problem of `trajectory planning' where 
the control solution is modeled as deviations from an underlying dynamical 
model~\citep{ijspeert2001trajectory,schaal2007dynamics}.
A dynamical systems view of optimally controlled trajectories can also be 
related to a free energy principle for the state cost~\citep{friston2010free}.
In this letter, we take a similar dynamical systems perspective and 
show that the immediate state cost function in the optimal control 
cost can be associated with a known least-action principle for a class of 
stochastic optimal control problems. 
The derived dynamical model immediately yields optimally controlled trajectories,
obviating the need for iterative solutions as in dynamic programming.

We consider control problems in the class of Linear Markov Decision 
Processes (LMDPs)~\citep{fleming1982optimal,Todorov11478,kappen2005linear}.
LMDPs get their name because the Hamilton-Jacobi-Bellman equation for the 
optimal cost-to-go becomes a linear differential operator after a suitable 
exponential transform of the cost-to-go function. In the continuous time case, 
LMDP models can be solved using path integral Monte Carlo 
techniques~\citep{pra1990markov,Kappen_2005,kappen2016adaptive}, while discrete time 
formulations can be solved using a variety of methods including eigenvalue 
problems and temporal difference reinforcement 
learning~\citep{todorov2006linearly}.
In this work, we exploit the relationship of LMDP models to the \schr 
equation~\citep{nagasawa1989transformations,pra1990markov,dai1991stochastic,pavon1991free,Kappen_2005,OHSUMI2019181} to 
reinterpret the control cost function in terms of a dynamic potential. 
Similar approaches linking integrable systems to stochastic optimal control of 
mean field games was presented in~\citet{swiecicki2016schrodinger} and control of ensembles of non-interacting entities in~\citet{bakshi2020schrodinger}.
% although here we do not assume Langevin dynamics as is common in these models. 
One potential advantage of introducing a \schr representation of 
optimal control is that the \schr solutions automatically explore all 
possible control paths, which can be seen in the path integral formulation 
of \schr equation solutions~\cite{chapline2001,Kappen_2005}

Quantum inverse scattering defines a one-to-one relationship between 
asymptotic `scattering data' and the potential~\citep[e.g.,][]{newton}.
\citet{rose1996global,rose2001single,Rose_2003} showed that Schr\"{o}dinger 
scattering solutions can be focused such that the transmitted scattering wave 
is a Dirac $\delta$-function at a specified later time. These 
`single-sided focusing' methods give a scattering interpretation to the 
topic of \schr bridges~\citep{chen2016relation}, which derive from a problem 
originally introduced by \schr~\citep{schrodinger1931umkehrung} to provide a 
probabilistic derivation of his wave equation.
\citet{carroll1990ubiquitous} derived a similar 
result by considering the relationships between spectral representations 
of related families of second-order differential operators.
\citet{Dyson:75,Dyson:1976cu} discovered 
a relationship between quantum inverse scattering methods in two-dimensions and
optimal feedback control of a noisy temporal signal.
The single-sided focusing method of Rose considered a known potential and sought
solutions for the incident wave that would be focused. 
In contrast, here we consider a known
incident wave and seek the potential, and thus implicitly the state cost, 
that causes the transmitted wave to be focused to a narrow peak. 

% ---------------------------------------------------------
\paragraph{\schr control problems}
\label{sec:schrodinger_control}

LMDP control problems are described by the following process model for 
a $p$-dimensional state $\xv$,
\begin{equation}\label{eq:system_dyn}
    d\xv = \bv(t, \xv(t))\, dt + \mathsf{C}\mathbf{a}(t, \xv(t))\, dt + d\xi,
\end{equation}
where $\mathbf{b}$ is a $p$-dimensional dynamical drift, 
$\mathsf{C}$ is a $p\times n$ control projection matrix,
$\mathbf{a}$ is an $n$-dimensional control variable, and 
$d\xi$ describes a $p$-dimensional Wiener process 
with $\mathbb{E}(d\xi d\xi) = \nu dt$ for a $p\times p$ covariance $\nu$.
We seek to find a control $\mathbf{a}(t,x(t))$, $t_0=0 < t < t_f$ such that the 
following cost function is mimized~\citep{Kappen_2005},
\begin{multline}\label{eq:control_cost}
    L(t_0, \xv_f, \mathbf{a}(\cdot)) \equiv 
    \mathbb{E}_{P(t_0,\xv_0)} \Biggl[
        q_f(\xv(t_f)) + \\
        \int_{t_0}^{t_f} dt\, \left(\frac{1}{2} \mathbf{a}^{\top}(t,\xv(t)) 
        	\mathsf{m}^{-1} \mathbf{a}(t,\xv(t))
        + q(t,\xv(t))\right)
    \Biggr],
\end{multline}
where $q$ is an arbitrary state cost function, $q_f$ is an asserted final time 
state cost, $\mathsf{m}$ is a matrix of control cost weights, and $P(t_0,\xv_0)$ 
denotes the distribution of paths under the dynamics of 
Equation~(\ref{eq:system_dyn}) that start at state $\xv_0$ at time $t_0$.
The action cost is a quadratic function of $\xv$, which can be interpreted 
as the continuous state space limit of a Kullback-Liebler divergence action 
cost that appears in the discrete state space formulation of 
LMDPs~\cite{todorov2008general}.

Minimizing the expected cost, $L$, with respect to the choice of actions 
defines the optimal cost-to-go function, which is the primary objective of the 
optimal control calculation,
\begin{equation}
 	J(t, \xv) \equiv \underset{\mathbf{a}(t\rightarrow t_f)}{\rm min} \
 	L(t, \xv, \mathbf{a}).
\end{equation} 
Then, substituting the optimal control variable $\mathbf{a}(t,\xv(t))$ leads to 
the Hamilton-Jacobi-Bellman (HJB) equation~\citep{kappen2005linear} that can be
used to determine the optimal cost-to-go function,
\begin{equation}\label{eq:HJB_opt}
    -\partial_t J = -\frac{1}{2}
    	\left(\nabla S\right)^{\top} \mathsf{C}\mathsf{m}^{-1}\mathsf{C}^{\top}
    	\nabla S
    + \bv\cdot\nabla J
    + {\rm Tr}\left(\frac{\nu}{2}\Delta J\right) + q,
\end{equation}
with boundary condition $J(t_f,\xv) = q_f(\xv)$.
By applying a Cole-Hopf transform~\citep{Hopf1950ThePD,Cole1951OnAQ}, we obtain 
an expression for the desirability function~\citep{Todorov11478}, which 
represents a partition function over optimally controlled stochastic 
paths~\citep{Kappen_2005},
\begin{equation}\label{eq:cole_hopf}
    z(t,\xv) \equiv e^{-J(t,\xv) / \lambda}.
\end{equation}
With this transform, the HJB equation becomes linear in $z$,
\begin{equation}\label{eq:linear_HJB}
    \partial_t z + \frac{1}{2}{\rm Tr}\left(\nu\Delta z\right) +  \mathbf{b} \cdot \nabla z 
        - \frac{q}{\lambda}z = 0.
\end{equation}
Here we have taken the noise covariance
$\nu = \lambda \mathsf{C}\mathsf{m}^{-1}\mathsf{C}^{\top}$, which 
determines the scalar parameter $\lambda$ in equation~(\ref{eq:cole_hopf}) 
from $\nu$, the control cost weights $\mathsf{m}$ and control projection 
matrix $\mathsf{C}$~\citep{Kappen_2005}.

As the next step in identifying the \schr equation related to this control problem,
define a new transform~\cite{nagasawa1989transformations,chapline2001,chapline2004quantum},
\begin{equation}\label{eq:cole_hopf_with_density}
	z(t,\xv) = \exp\left(R(t,\xv) - S(t,\xv)/\lambda\right).
\end{equation}
Substituting Equation~(\ref{eq:cole_hopf_with_density}) back into 
Equation~(\ref{eq:linear_HJB}) we find the sum of two expressions,
\begin{multline}\label{eq:linear_HJB_with_density}
	0 = -\frac{1}{\lambda}\Biggl[
		\partial_t S + \bv \cdot \nabla S - 
		\frac{1}{2}\left(\nabla S\right)^{\top} 
			\mathsf{C}\mathsf{m}^{-1}\mathsf{C}^{\top} \nabla S \\
		+ \frac{1}{2}{\rm Tr}\left(\nu\Delta S\right) + \tilde{q}
	\Biggr] \\
	+ \Biggl[\partial_t R + (\bv - \mathbf{v})\cdot\nabla R\Biggr],
\end{multline}
where $\tilde{q} \equiv q + \mathcal{B}$,
\begin{equation}
	\mathcal{B} \equiv -\frac{\lambda^2}{2}\left[
		{\rm Tr}\left(\mathsf{m}^{-1}\Delta R\right) 
		+ \left(\nabla R\right)^{\top}\mathsf{m}^{-1}\nabla R
	\right],
\end{equation}
is familiar as the Bohm potential in quantum 
mechanics~\citep{bohm1952suggested}, and 
$\mathbf{v} \equiv \mathsf{m}^{-1}\nabla S$ is a current velocity, familiar from 
Nelson stochastic mechanics~\citep{nelson2020dynamical}.
The expression in the first set of brackets of 
Equation~(\ref{eq:linear_HJB_with_density}) 
is thus one side of an HJB equation for $S$ but with a cost function that is 
modified from that in the HJB equation for $J$ by the addition of the Bohm potential.
The expression in the second set of brackets of 
Equation~(\ref{eq:linear_HJB_with_density})
is one side of a continuity equation, which we interpret in terms of a state 
probability density $\mu \equiv \exp(2R)$.
Equation~(\ref{eq:linear_HJB_with_density}) is satisfied if the expressions in 
each line are separately equal to zero, giving an HJB equation for $S$ and a 
continuity equation for $\mu$. But, in general we only require the sum to be 
zero so that lack of conservation of state probability can be compensated for 
by proportional deviations from HJB optimality for the cost-to-go function $S$.

Because the Bohm potential is equal to the Fisher information for the density 
$\mu$, we can see that high Fisher information, corresponding to a narrow 
density, incurs more control cost than low information, corresponding to a 
broad density. This trend is consistent with the theory of the linear Bellman 
equation~\citep{Kappen_2005} where control becomes more expensive in regions 
of state space less likely to be visited by the process noise. 

With Equation~(\ref{eq:cole_hopf_with_density}) and separation of the two lines of 
Equation~(\ref{eq:linear_HJB_with_density}), we have arrived at an 
interpretation of the desirability function $z$ in terms of a diffusion process 
with state density $\mu$ that is controlled with optimal cost-to-go $S$.
\citet{nagasawa1989transformations} shows that just such a diffusion process can 
be equivalently described by the \schr equation upon making the identifications,
\begin{align}
	z = \exp(R - S/\lambda) &\leftrightarrow \psi = \exp(R - iS/\lambda), \\
	\hat{z} = \exp(R + S/\lambda) &\leftrightarrow \psi^{\dagger} = \exp(R + iS/\lambda), \notag \\	
	z\hat{z} &= \mu = \psi\psi^{\dagger}.
\end{align}
The wavefunction $\psi$ thus defined satisfies
the \schr equation for a particle in a magnetic field,
\begin{multline}\label{eq:schrodinger_eq}
	i\hbar\partial_t \psi(t,\xv) = 
	-\frac{\hbar}{2}{\rm Tr}\left(\nu \Delta \psi\right) \\
	+ i\hbar \bv(t,\xv)\cdot\nabla \psi
	+ V(t,\xv)\psi,
\end{multline}
where $\bv$ is now interpreted as a magnetic vector potential and $V$ is a new 
scalar potential that is related to the state cost 
$\tilde{q}$~\citep[][equation 4.14]{nagasawa1989transformations},
\begin{equation}\label{eq:potential_from_cost}
	V = \tilde{q} - 2\partial_t S - \nu (\nabla S)^2 - 2\bv\cdot\nabla S.
\end{equation}
By solving the \schr equation with this potential $V$, we obtain both functions 
$R$ and $S$, which yield the solution to the optimal control problem in the 
class of linear Markov decision processes solved by the linear Bellman equation.

% ---------------------------------------------------------
\paragraph{Focusing principle}
\label{sec:focusing_principle}

The relation in Equation~(\ref{eq:potential_from_cost}), if interpreted directly,
makes the \schr equation nonlinear and thus difficult to solve. 
However, if the potential $V$ is known, then we can more easily solve the 
linear \schr equation.

The theory of quantum inverse scattering describes the problem of determining 
an unknown potential from spatially asymptotic information about the 
wavefunction, called the `scattering data'~\cite{morse1954methods}. 
We will identify the scattering 
data with the starting and final target state distributions for 
a finite horizon optimal control problem.

The unique potential that corresponds to given scattering data can be determined 
from a least-action principle where the action is defined as the integrated 
squared amplitude of the `tail' of the transmitted incident 
wave~\citep{rose1996global,Rose_2003},
\begin{equation}\label{eq:blurring_action}
	A \equiv \int d\xv\, \left[\phi(t_f,\xv) - \phi_{\rm in}(\xv)\right]^2,
\end{equation}
where $\phi_{\rm in}$ denotes an incident wavefront and $\phi(t,\xv)$ is the 
transmitted wave.
By minimizing this action, an incident wave is focused 
to a Dirac $\delta$-function upon passing over the potential, which is called 
single-sided focusing~\citep{rose2001single}.

The Euler-Lagrange equation for the action in Equation~(\ref{eq:blurring_action})
is the Mar\v{c}enko integral equation~\citep{Rose_2003}, which in one dimension 
is,
\begin{equation}
    \Omega(-\tau;x_f) + \Gamma(\tau+x_f) + \int_{-\infty}^{\infty} d\tau'\,
        \Gamma(\tau + \tau') \Omega(-\tau'; x_f) = 0,
\end{equation}
for $\tau < x_f$ and $\Omega(-\tau)=0$ for $\tau > x_f$ and where 
$\Gamma$ is the inverse Fourier transform of the reflection 
coefficient (assuming no bound states) and $\Omega$ is a causal kernel that 
is related to the potential as,
\begin{equation}
	V(x) = 2\partial_x \Omega(x,x).
\end{equation}
Focusing in two-dimensions can be derived from an analogous two-dimensional Marcenko equation~\citep{,cheney1984inverse,cheney1985two,Yagle_1998}.

In \citep{rose2001single}, it is shown that 
maximum focusing corresponds to focusing 
of the real part of the wavefunction $\psi$.
Focusing the real part of the wavefunction implies focusing also of the 
probability current density, $j$, so that,
\begin{equation}
    j(t_f,x) \equiv e^{2R(t_f,x)} \partial_x S(t_f,x) \approx m\delta(x - x_f).
\end{equation}
Because the cost-to-go $S(t_f,x)$ at the final time is asserted as a boundary 
condition, the function $\partial_x S$ can be arbitrary. The only way 
to ensure the current density $j$ is focused is if the state probability 
density is focused,
% From the nature of our boundary conditions, we have,
\begin{equation}
    \mu(t_f,x) = e^{2R(t_f,x)} = \delta(x - x_f).
\end{equation}

By using the focusing principle to derive $V$ from given boundary conditions 
for $\psi$, we can obtain $\tilde{q}$ from Equation~(\ref{eq:potential_from_cost})
and then get the original state cost $q = \tilde{q} - \mathcal{B}$.
Thus, for the stochastic finite-horizon optimal control problem defined 
in Equations~(\ref{eq:system_dyn}) and (\ref{eq:control_cost})
with convex state cost $q_f(x)$ at final time $t_f$,
there exists a unique $q(x,t)$ that minimizes the variance in the 
final state at $t=t_f$ relative to an asserted target value at $x=x_f$.

% ---------------------------------------------------------
\paragraph{Numerical algorithm}
\label{sec:algorithm}

Inspired by the action in Equation~(\ref{eq:blurring_action}), 
we define a `focusing metric' for numerical optimization of 
the potential $V$,
\begin{equation}\label{eq:focusing_metric}
	\mathcal{F}(V) \equiv \sum_{x} \left[
		{\rm Re}\left(\psi(t_f,\xv; V)\right) - {\rm target}(\xv)
	\right]^2,
\end{equation}
where $\psi$ is required to be a solution of the \schr 
equation~(\ref{eq:schrodinger_eq}) and the `target' is an asserted function that is square-normalizable. 

The optimization of the metric $\mathcal{F}$ with respect to $V$ can 
be accomplished by solving the inverse scattering problem via
numerical integration of the Mar\v{c}enko equation. However, obtaining the 
reflection coefficient (and any bound states) from the initial state distribution
can be numerically challenging. Instead, we use gradient descent to optimize 
the metric $\mathcal{F}$, where the gradient operates through the numerical 
eigensolver for the time-independent \schr equation using the Jax 
software library~\citep{jax2018github}.

% ---------------------------------------------------------
\paragraph{Numerical example}
\label{sec:examples}

To demonstrate the focusing algorithm, we adopt the common test problem of an agent 
navigating a simple maze in two dimensions. The maze is composed of both 
exterior walls and interior 
barriers. The agent is positioned in the top left corner of the maze at $t_0=0$
and must reach the bottom right corner of the maze by $t_f=0.6$. The agent path 
is terminated if the agent touches a wall or barrier.

We solve the control problem on an arbitrary small test system consisting of a
discretized spatial grid of $51\times51$ cells
over an extent $-1 \le x_1, x_2 \le 1$. We take $\hbar = \lambda = 1$ and 
$m = 0.5$.
We solve the \schr equation by first finding the energy eigenfunctions and 
eigenvalues for a given potential $V$. We then expand the wavefunction in 
eigenfunctions, truncating to the 15 eigenfunctions with the smallest eigenvalues.
We initialized $V$ to a quadratic function of distance on the grid measured 
from the starting location in the top left. We then performed gradient descent
with a learning rate of 0.02 until the slope of the learning curve flattened out.

\begin{figure}[!htb]
	\centerline{
		\includegraphics[width=0.5\textwidth]{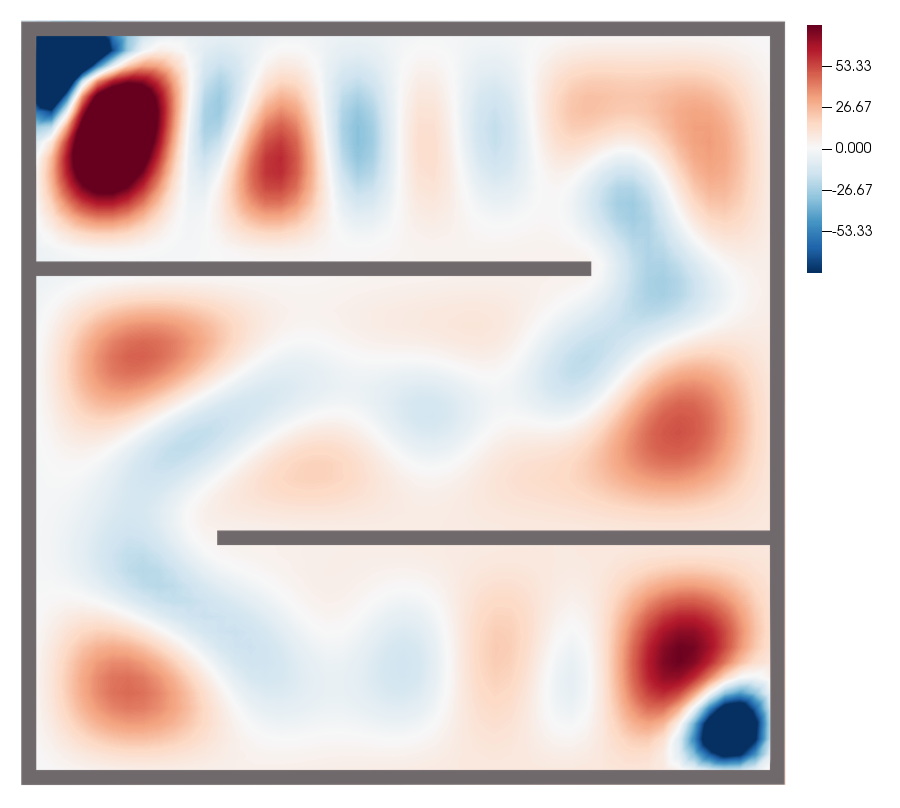}
	}
	\caption{\label{fig:gridworld_potential}
		The potential that approximately minimizes the focusing metric for a two-dimensional 
		gridworld maze. The agent begins in the top left of the grid and must navigate to the 
		bottom right without hitting any of the exterior or interior walls.
	}
\end{figure}
\begin{figure}[!htb]
	\centerline{
		\hspace*{-0.2cm}
		\includegraphics[width=0.147\textwidth]{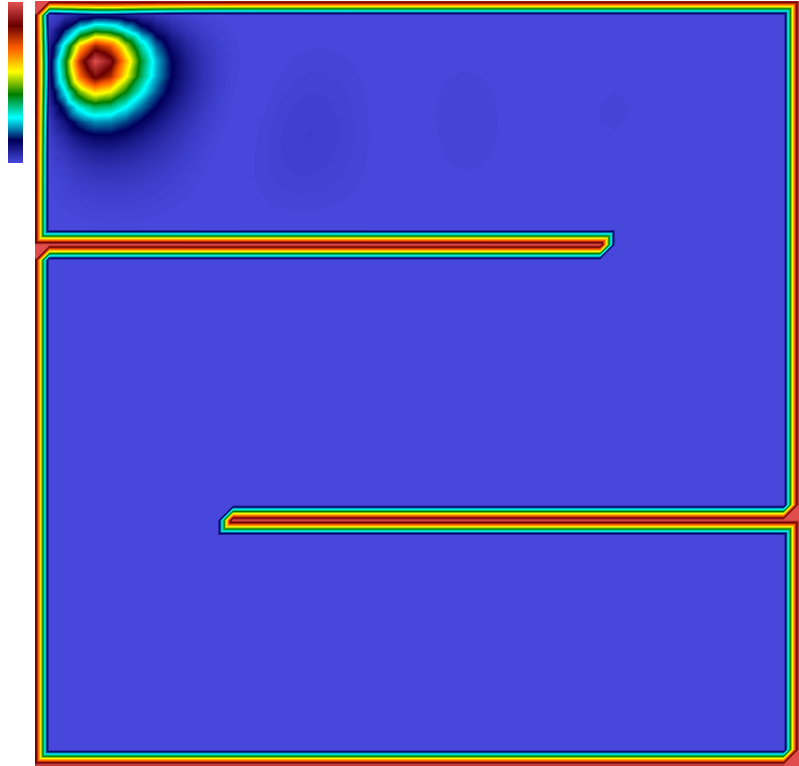}
		\includegraphics[width=0.14\textwidth]{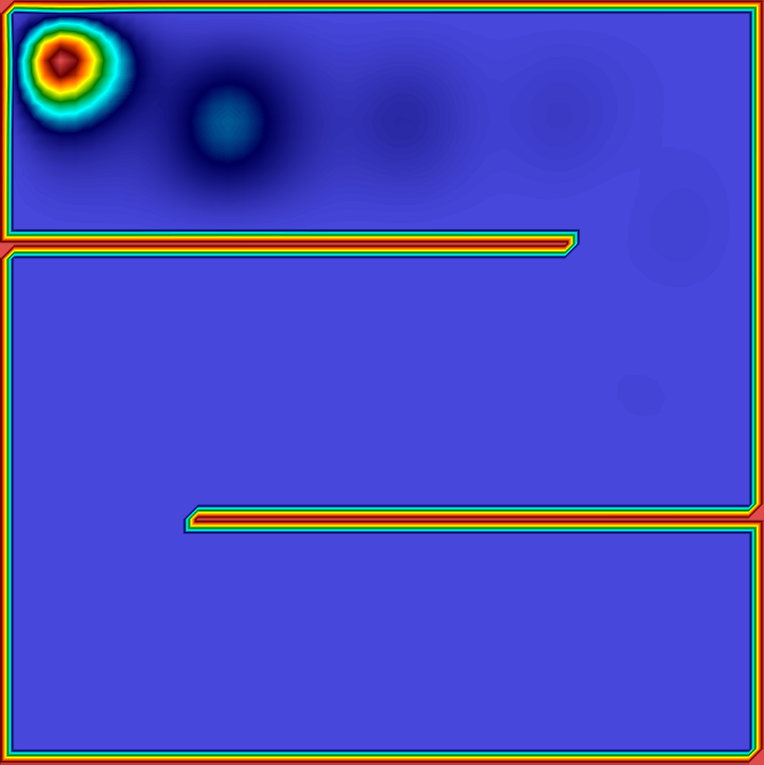}
		\includegraphics[width=0.14\textwidth]{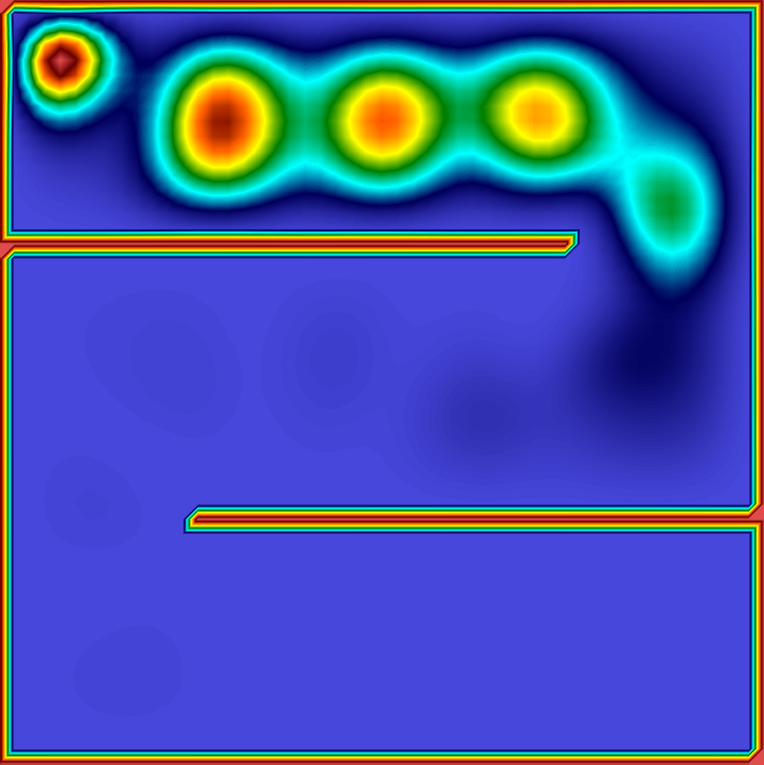}
	}
	\centerline{
		\includegraphics[width=0.14\textwidth]{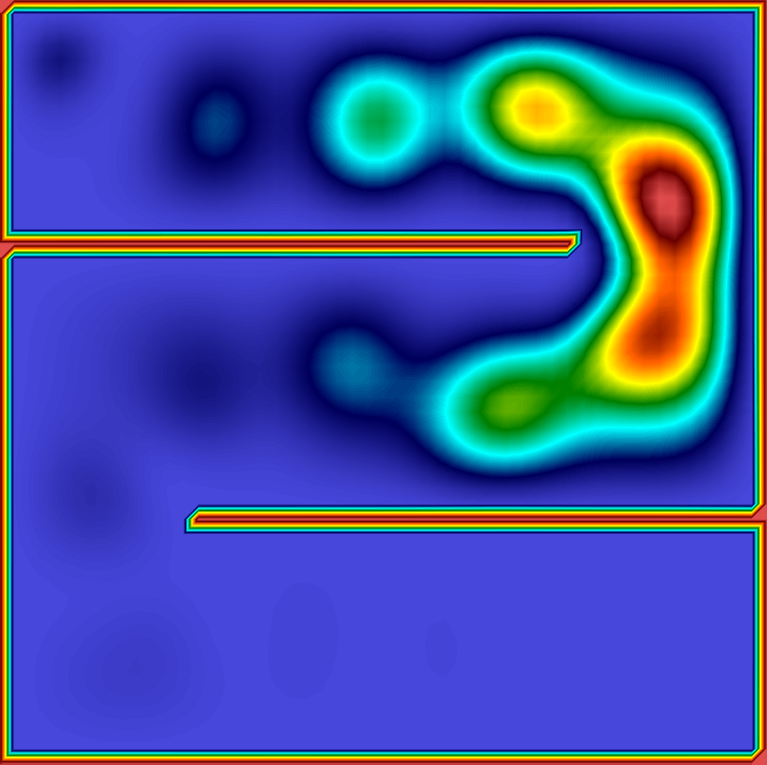}
		\includegraphics[width=0.14\textwidth]{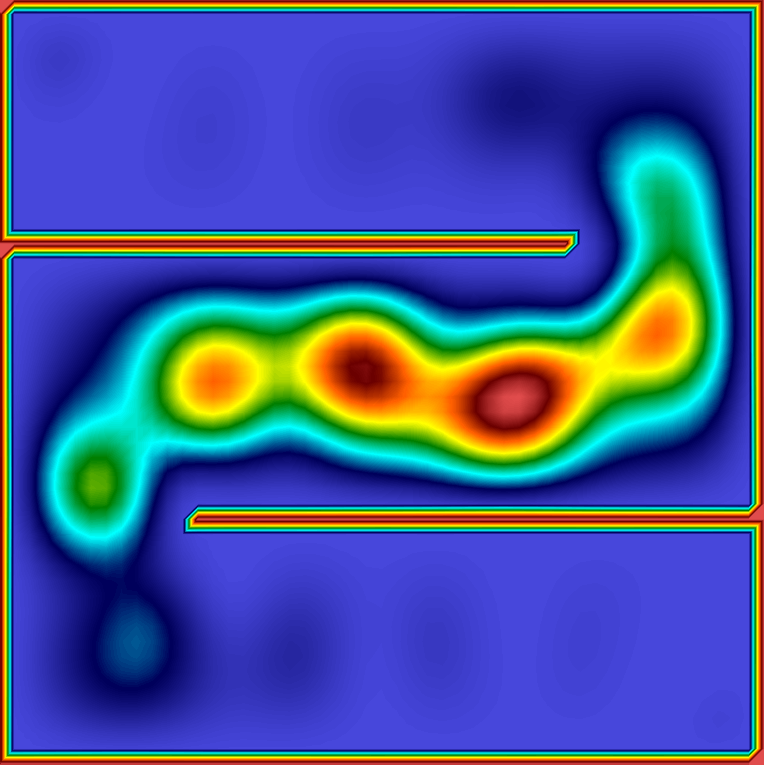}
		\includegraphics[width=0.14\textwidth]{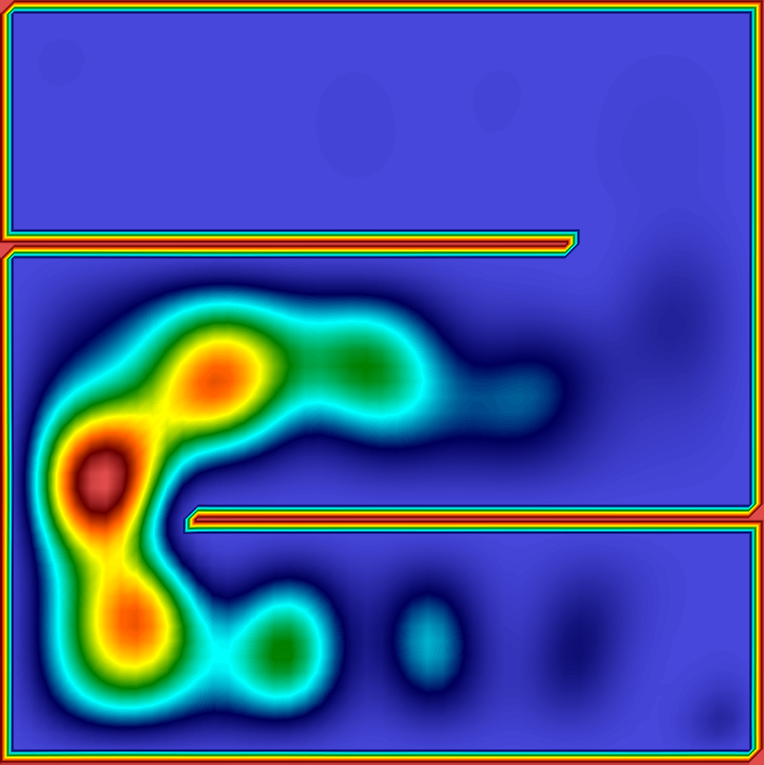}
	}
	\centerline{
		\includegraphics[width=0.14\textwidth]{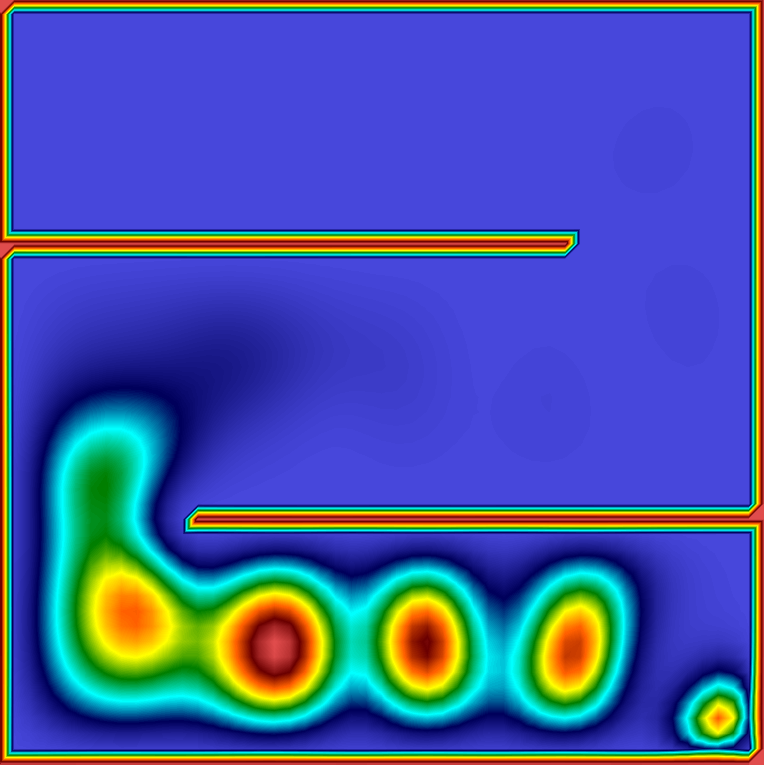}
		\includegraphics[width=0.14\textwidth]{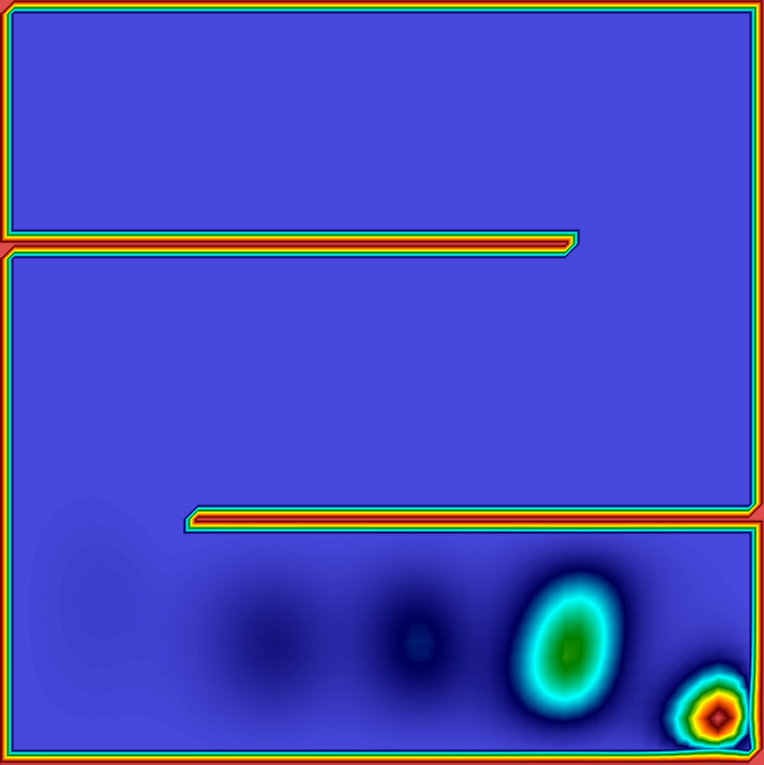}
		\includegraphics[width=0.14\textwidth]{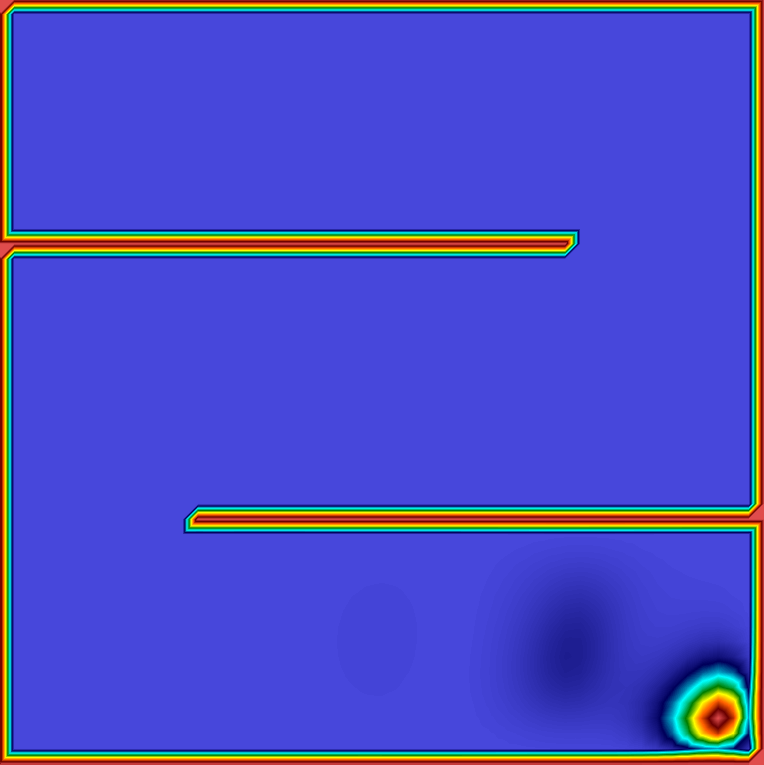}
	}
	\caption{\label{fig:gridworld_density}
		The agent state probability density as a function of time for the 
		optimally controlled solution to navigating a gridworld maze. Time evolves in each panel from top left to bottom right.
	}
\end{figure}

We show the resulting potential $V$ that approximately minimizes $\mathcal{F}$ 
in Figure~(\ref{fig:gridworld_potential}) and the state probability 
distribution over time for the optimally-trained agent in 
Figure~(\ref{fig:gridworld_density}). We see that the probability density 
focuses to a narrow target distribution at the final time as desired. In 
addition, the probablity density at all times navigates the barriers effectively
to ensure zero probability of early termination of the agent maze navigation
episode.

% ---------------------------------------------------------
\paragraph{Conclusions}
\label{sec:conclusions}

We have shown that we can solve stochastic optimal control problems through computation of a potential function with an inverse scattering method.
We assume we are given knowledge about the environment at the current time step and an asserted final time goal that the probability distribution of the agent state is `focused' on specified state space locations.
The control problem is then solved by computing the potential for all states and times that serves to drive an agent towards the goal. The potential is computed by matching the wavefunction solution of the \schr equation to a desired target, which is a computation that can be executed in parallel across states and times. 

The agent in this formulation is viewed as a dynamical system that follows `forces' equal to the gradients of the potential function. Equivalently, the HJB optimal cost-to-go for the agent can be obtained from the phase of the wavefunction. While not demonstrated here, the approach admits a trivial generalization to multi-agent systems, where each agent learns its own potential.

% We see from Figure~\ref{fig:gridworld_potential} that the potential is localized in state space ({\it i.e.}, has a finite correlation length) as expected for a solution of the GLM equation. This localization means that control actions in future implementations might be computed in parallel across state and time locations. 

Single-sided focusing minimizes the path length (in phase units) to reach a desired end state. Various paths considered in the variational context define a surface embedded in state space, whose area is the loss of information obtained by the agent relative to the optimal path~\citep{chapline2001}. For small deviations from the optimal path, the state cost function is determined by the probability distribution of the agent state, giving a localized description of the entire control solution that may be easily parallelized in future computational applications and can be contrasted with non-local Euclidean path integral methods~\citep{kappen2016adaptive} or episode-based training in RL~\citep{sutton2018reinforcement}. 

The focusing action in Equation~(\ref{eq:blurring_action}) is only one possible
choice for determining the potential. Indeed, other actions are known to produce
Euler-Lagrange equations that describe integrable systems, such as the 
Korteweg-de Vries (KdV) equation~\citep{Zakharov1971}. The optimal cost 
functions for control problems might then be derived as solutions to a broader 
class of integrable models. 
In the most general case of a time-dependent potential, the \schr 
equation can be interpreted as one of the Lax operators~\cite{lax1968integrals} 
for a completely integrable dynamical model for the potential. This is another 
way that integrability can appear as a defining feature for solutions of 
stochastic optimal control problems.

% ---------------------------------------------------------
\begin{acknowledgments}
This work was performed under the auspices of the U.S. Department of Energy 
by Lawrence Livermore National Laboratory under Contract DE-AC52-07NA27344.
Funding for this work was provided by LLNL Laboratory Directed Research and 
Development grant 22-SI-001.
\end{acknowledgments}

% ---------------------------------------------------------
% \bibliography{ds}
% \input{paper_prl.bbl}
% ---------------------------------------------------------


\end{document}